\documentclass[]{article}

\usepackage{amsmath,latexsym,amssymb}
\usepackage{listings}
\newtheorem{theorem}{Theorem}
\newtheorem{definition}{Definition}

\DeclareMathOperator{\erfc}{erfc}
\title{On Computation of Matrix\\Mittag-Leffler Function}
\author{Ivan Matychyn}
\date{}

\begin{document}

\maketitle

\begin{abstract}
A method for computation of the matrix Mittag-Leffler function is presented. The method is based on Jordan canonical form and implemented as a Matlab routine \cite{matlab}.
\end{abstract}

\section{Fractional Differential Equations and Matrix Mittag-Leffler Functions}

The matrix Mittag-Leffler function was probably first introduced in the paper \cite{ChikEid}, where it was used in an explicit solution of a linear system of fractional order equation (FDEs)
\begin{equation}\label{eq1}
D^\alpha z=Az + f,\ 0<\alpha \le 1.
\end{equation}

Here $ D^\alpha z$ stands for the Riemann--Liouville fractional derivative of order $ \alpha$. In general, if $ g$ is a function having absolutely continuous derivatives up to the order $ m-1$, the Riemann--Liouville derivative of fractional order $ \alpha$, $m-1<\alpha\leq m$, can be defined as follows:
\begin{equation}\label{rl}
D^\alpha g(t)=\frac{1}{\Gamma(m-\alpha)}\frac{d^m}{dt^m} \int_0^t \frac{g(\tau)}{(t-\tau)^{\alpha-m+1}}d\tau.
\end{equation}

Hereafter $ A$ is a fixed real $ n\times n$ matrix, and $ z, f:\ [0,\infty)\to \mathbb{R}^n$ are measurable vector-functions taking values in $ \mathbb{R}^n$. 

If \eqref{eq1} is supplied with initial condition of the form
\begin{equation}\label{eq2}
\frac{1}{\Gamma(1-\alpha)}\int_{0}^{t} \frac{z(\tau)}{(t-\tau)^\alpha} d\tau \Bigr|_{t=0}=z^0,
\end{equation}
then solution to the initial value problem \eqref{eq1}, \eqref{eq2}  can be written down in the form
\begin{equation}\label{eq3}
z(t)=t^{\alpha-1} E_{\alpha,\alpha}(At^\alpha)z^0 + \int_{0}^{t}(t-\tau)^{\alpha-1} E_{\alpha,\alpha}(A(t-\tau)^\alpha) f(\tau) d\tau,
\end{equation}
where 
\begin{equation} \label{mlf} 
E_{\alpha,\beta}(A)=\sum_{k=0}^{\infty}\frac{A^k}{\Gamma(\alpha k + \beta)},\ \alpha>0,\ \beta\in\mathbb{C},
\end{equation}
denotes the matrix Mittag-Leffler function of $A$.

The  expression \eqref{eq3} can be rewritten in more compact form
\begin{equation}\label{eq4}
z(t)=e_\alpha^{At}z^0 + \int_{0}^{t} e_\alpha^{A(t-\tau)} f(\tau) d\tau,
\end{equation}
where $ e_\alpha^{At}=t^{\alpha-1} E_{\alpha,\alpha}(At^\alpha)$ is the matrix $\alpha$-exponential function introduced in the monograph \cite{KilSriTru}.

Since FDEs involving the Riemann--Liouville fractional derivative require initial conditions of the form \eqref{eq2} lacking clear physical interpretation, the regularized fractional derivative was introduced. The latter is often referred to as the Caputo derivative and defined as follows:
\begin{equation}\label{caputo}
D^{(\alpha)}g(t)=\frac{1}{\Gamma(m-\alpha)} \int_0^t \frac{g^{(m)}(\tau)}{(t-\tau)^{\alpha-m+1}}d\tau,\ m-1<\alpha\leq m.
\end{equation}

Initial value problem for FDEs involving the Caputo derivative
\begin{equation}\label{fdecaputo}
D^{(\alpha)} z=Az + f,\ 0<\alpha \le 1,
\end{equation}
requires standard initial conditions
\begin{equation}\label{inicaputo}
z(0)=z^0,
\end{equation}
and its solution can be explicitly written down in terms of matrix Mittag-Leffler functions as follows \cite{ChikMat}:
\begin{equation}\label{solcaputo}
z(t)=E_{\alpha,1}(At^\alpha)z^0 + \int_{0}^{t}(t-\tau)^{\alpha-1} E_{\alpha,\alpha}(A(t-\tau)^\alpha) f(\tau) d\tau.
\end{equation}

As an example let us consider well-known Bagley--Torvik equation \cite{BagTor} describing vibrations of a rigid plate immersed in Newtonian liquid:
\begin{equation} \label{bt1}
ay''(t) + b D^{(3/2)}y(t) +cy(t)=f(t)
\end{equation}
\begin{equation} \label{bt2} 
y(0)=y_0,\quad{} y'(0)=y'_0.
\end{equation}
Its analytical solution obtained with the help of fractional Green's function in terms of scalar generalized Mittag-Leffler functions is cumbersome and involves evaluation of a convolution integral, containing a Green's function expressed as an infinite sum of derivatives of Mittag-Leffler functions, and for general functions $f$ this cannot be evaluated conveniently.

The equation of Bagley--Torvik is equivalent to the following system \cite{Diet}
\begin{equation}
\label{eq:sys}
D^{(1/2)}z=Bz+Cf,
\end{equation}
where
$B=\left[ \begin{array}{cccc} 
0 &1 &0 &0\\
0 &0 &1 &0\\
0 &0 &0 &1\\
-c/a &0 &0 &-b/a\\
\end{array} \right]$,
$C=\left[ \begin{array}{c} 0\\ 0\\ 0\\ 1/a \end{array} \right]$, $z=\text{col}(y, D^{(1/2)}y, y', D^{(3/2)}y )$ under the initial conditions
\begin{equation*}
z(0)=z_0=\text{col}( y_0, 0, y'_0, 0).
\end{equation*}
Its solution in terms of matrix Mittag-Leffler functions is given by the following expression:
\begin{equation} \label{bt4}
z(t)=E_{\frac{1}{2},1}(B\sqrt{t})z_0+\int_0^t  E_{\frac{1}{2},\frac{1}{2}}\left(B\sqrt{t-\tau}\right) C\frac{f(\tau)d\tau}{\sqrt{t-\tau}},
\end{equation}
which can be easily evaluated.

The explicit expressions \eqref{eq3}, \eqref{eq4}, and \eqref{solcaputo} play a key role in numerous applications related to systems with fractional dynamics \cite{app1,app2,app3}. That is why the methods for computing the matrix Mittag-Leffler function are so important.

Both the matrix Mittag-Leffler function and the matrix $\alpha$-exponential functions are generalizations of matrix exponential function, since
\begin{equation*}
E_{1,1}(At)=e_1^{At}=e^{At}.
\end{equation*}
This implies that some of numerous existing methods for computing the matrix exponential can be adapted for the matrix Mittag-Leffler functions as well. An overview and analysis of these methods can be found in the paper \cite{MolerVanLoan} and in the monograph \cite{Higham}. Unfortunately, the technique of scaling and squaring, widely used in computing of the matrix exponential, cannot be applied for the matrix Mittag-Leffler and $\alpha$-exponential functions, as the latter do not possess the semigroup property.

Here we describe a method of computing the matrix Mittag-Leffler function based on  the Jordan canonical form representation. This method is implemented with \textsc{Matlab} code \cite{matlab}.

\section{Matrix Functions}
There exists a number of equivalent definitions of a matrix function. 
The following classic definition in terms of interpolation polynomials is according to \cite{gantmacher1959theory}.  Let
\begin{equation*}
\psi(\lambda)=(\lambda-\lambda_1)^{m_1}(\lambda-\lambda_2)^{m_2}\ldots(\lambda-\lambda_s)^{m_s}
\end{equation*}
be the minimal polynomial of $ A$, where $\lambda_1, \lambda_2,  \ldots, \lambda_s$  are all the distinct eigenvalues of $ A$. The degree of this polynomial is $ m=\sum_{k=1}^{s} m_k$.

Let us consider a sufficiently smooth function $ f(\lambda)$ of scalar argument and call the $ m $ numbers
\begin{equation} \label{eq5}
f(\lambda_k), f'(\lambda_k), \ldots, f^{(m_k-1)}(\lambda_k) \quad (k=1,\ldots,s)
\end{equation} 
the values of the function $ f$ on the spectrum of the matrix $ A$ and the set of all these values will be denoted
symbolically by $ f(\Lambda_A)$. If for some function $ f $ the values \eqref{eq5} exist, then we will say that the function $ f $ is defined on the spectrum of the matrix $ A $.
\begin{definition}[matrix function via interpolation polynomial \cite{gantmacher1959theory}] \label{inrpoldef}
Let $ f(\lambda) $ be a function defined on the spectrum of a matrix $ A$ and $ r(\lambda) $ the corresponding interpolation polynomial such that $ f(\Lambda_A)=r(\Lambda_A)$. Then
\[
f(A)=r(A).
\]
\end{definition}

Let us recall the following well-known
\begin{theorem}
Any constant $ n\times n$ matrix $ A $ is similar to a matrix $ J $ in Jordan canonical form. That is, there exists an invertible matrix $ P $ such that the $ n\times n$ matrix $J= Z^{-1} A Z$ is in the canonical form 
\begin{equation} \label{jordan}
J=\mathrm{diag} \{J_1,J_2,\ldots,J_s\}
\end{equation}
where each Jordan block matrix $ J_k$, $ k=1,\ldots,s $, is a square matrix of the form
\[
J_k=
\begin{pmatrix}
\lambda_k & 1 & 0 & \ldots & 0\\
0 & \lambda_k & 1 & \ldots & 0\\
0 & 0 &\lambda_k & \ldots & 0\\
\vdots & \vdots & \vdots & \ddots & \vdots\\
0 & 0 & 0 & \ldots & 1\\
0 & 0 & 0 & \ldots & \lambda_k
\end{pmatrix}.
\]
\end{theorem}

It is shown (see e.g. \cite{gantmacher1959theory}) that Definition \ref{inrpoldef} is equivalent to the following definition based on the Jordan canonical form. The latter we will use for computing the matrix Mittag-Leffler function. 

\begin{definition}[matrix function via Jordan canonical form] \label{jordandef}
Let the function $ f $ be defined on the spectrum of $ A $ and let $ A=Z J Z^{-1} $, where $ J $ is the Jordan canonical form \eqref{jordan}. Then
\begin{equation}\label{eq6}
f(A)=Z f(J) Z^{-1}= Z\; \mathrm{diag} \{f(J_1),f(J_2),\ldots,f(J_s)\} Z^{-1},
\end{equation} 
where
\begin{equation}\label{eq7}
f(J_k)=
\begin{pmatrix}
f(\lambda_k) & f'(\lambda_k) & \frac{f''(\lambda_k)}{2} & \ldots & \frac{f^{(m_k-1)}(\lambda_k)}{(m_k-1)!}\\
0 & f(\lambda_k) & f'(\lambda_k) & \ldots & \frac{f^{(m_k-2)}(\lambda_k)}{(m_k-2)!}\\
0 & 0 & f(\lambda_k) & \ldots & \frac{f^{(m_k-3)}(\lambda_k)}{(m_k-3)!}\\
\vdots & \vdots & \vdots & \ddots & \vdots\\
0 & 0 & 0 & \ldots & f'(\lambda_k)\\
0 & 0 & 0 & \ldots & f(\lambda_k)
\end{pmatrix}.
\end{equation}
\end{definition}

\subsection{Generalized Mittag-Leffler Functions}
The generalized (scalar) Mittag-Leffler function also known as Prabhakar function is defined for complex $ z, \alpha, \beta, \rho \in \mathbb{C} $, and $\Re(\alpha) >0  $ by
\begin{equation} \label{gmlf} 
E_{\alpha,\beta}^\rho(z)=\sum_{k=0}^{\infty}\frac{(\rho)_k}{\Gamma(\alpha k + \beta)}\frac{z^k}{k!},
\end{equation}
where $(\rho)_k = \rho(\rho+1)\ldots (\rho+k-1)$ is the Pochhammer symbol.

In particular, when $ \rho = 1 $, it coincides with the Mittag-Leffler function \eqref{mlf}:
\[
E_{\alpha,\beta}^1 (z)=E_{\alpha,\beta} (z).
\]

Since the expression \eqref{eq7} involves derivatives, the following equation \cite{KilSriTru} is important for the purpose of computing the matrix Mittag-Leffler function:
\begin{equation}\label{mlfdrv}
\left( \frac{d}{dt} \right)^m E_{\alpha,\beta}(t)= m! E_{\alpha,\beta+\alpha m}^{m+1}(t), \quad m\in\mathbb{N}.
\end{equation}

In view of \eqref{mlfdrv}, the formulas \eqref{eq6}, \eqref{eq7} take on the form
\begin{equation}\label{qq6}
E_{\alpha,\beta}(A)= Z\; \mathrm{diag} \{E_{\alpha,\beta}(J_1),E_{\alpha,\beta}(J_2),\ldots,E_{\alpha,\beta}(J_s)\} Z^{-1},
\end{equation} 
\begin{equation}\label{qq7}
E_{\alpha,\beta}(J_k)=
\begin{pmatrix}
E_{\alpha,\beta}(\lambda_k) & E_{\alpha,\beta+\alpha}^2(\lambda_k) & E_{\alpha,\beta+2\alpha}^3(\lambda_k) & \ldots & E_{\alpha,\beta+(m_k-1)\alpha}^{m_k}(\lambda_k)\\
0 & E_{\alpha,\beta}(\lambda_k) & E_{\alpha,\beta+\alpha}^2(\lambda_k) & \ldots & E_{\alpha,\beta+(m_k-2)\alpha}^{m_k-1}(\lambda_k)\\
0 & 0 & E_{\alpha,\beta}(\lambda_k) & \ldots & E_{\alpha,\beta+(m_k-3)\alpha}^{m_k-2}(\lambda_k)\\
\vdots & \vdots & \vdots & \ddots & \vdots\\
0 & 0 & 0 & \ldots & E_{\alpha,\beta+\alpha}^2(\lambda_k)\\
0 & 0 & 0 & \ldots & E_{\alpha,\beta}(\lambda_k)
\end{pmatrix}.
\end{equation}

\section{Software Implementation}
The formulas \eqref{qq6}, \eqref{qq7} can be used for computing the matrix Mittag-Leffler function and were implemented in the form of \textsc{Matlab} routine \texttt{mlfm.m} \cite{matlab}. For computing of generalized Mittag-Leffler functions of the form $E_{\alpha,\beta+(m-1)\alpha}^{m}(\lambda_k)$, the \textsc{Matlab} routine by R.~Garrappa is used, which implements the optimal parabolic contour (OPC) algorithm described in \cite{Garrappa} and based on the inversion of the Laplace transform on a parabolic contour suitably chosen in one of the regions of analyticity of the Laplace transform.

To verify the accuracy of the \texttt{mlfm.m} routine, one can consider the Bagley--Torvik equation \eqref{bt1}, \eqref{bt2}. 

If $ c=0$, the matrix Mittag-Leffler functions appearing in \eqref{bt4} can be found analytically. 

Indeed, if $ c=0 $ the matrix $ B $ in \eqref{eq:sys} takes on the form 
\begin{equation} \label{B}
B=\left[ \begin{array}{cccc} 
0 &1 &0 &0\\
0 &0 &1 &0\\
0 &0 &0 &1\\
0 &0 &0 &p\\
\end{array} \right],
\end{equation}
where $p= -b/a$.

Hence,
\[
B^2=\begin{pmatrix} 
0 &0 &1 &0\\
0 &0 &0 &1\\
0 &0 &0 &p\\
0 &0 &0 &p^2\\
\end{pmatrix},\quad
B^k=\begin{pmatrix} 
0 &0 &0 &p^{k-3}\\
0 &0 &0 &p^{k-2}\\
0 &0 &0 &p^{k-1}\\
0 &0 &0 &p^k\\
\end{pmatrix},\quad k=3,4,\dots
\]
Therefore,
\begin{gather*}
E_{\frac{1}{2},1}(B)  = \sum\limits_{k=0}^\infty \frac{B^k}{\Gamma(k/2 + 1)}=\\ 
= 
\begin{pmatrix} 
1 &\frac{1}{\Gamma(3/2)} &1 &p^{-3}E_{\frac{1}{2},1}(p) - p^{-3} - \frac{1}{p^2\Gamma(3/2)} - \frac{1}{p}\\
0 &1 &\frac{1}{\Gamma(3/2)} &p^{-2}E_{\frac{1}{2},1}(p) - p^{-2} - \frac{1}{p\Gamma(3/2)}\\
0 &0 &1 &p^{-1}E_{\frac{1}{2},1}(p) - p^{-1}\\
0 &0 &0 &E_{\frac{1}{2},1}(p)\\
\end{pmatrix},\\
\\
E_{\frac{1}{2},\frac{1}{2}}\left(B\right) = \sum\limits_{k=0}^\infty\frac{B^k}{\Gamma((k+1)/2)}=\\ 
= 
\begin{pmatrix} 
\frac{1}{\Gamma(1/2)} &1 &\frac{1}{\Gamma(3/2)} &p^{-3}E_{\frac{1}{2},\frac{1}{2}}(p) - \frac{1}{p^3\Gamma(1/2)} - \frac{1}{p^2} - \frac{1}{p\Gamma(3/2)}\\
0 &\frac{1}{\Gamma(1/2)} &1 &p^{-2}E_{\frac{1}{2},\frac{1}{2}}(p) - \frac{1}{p^2\Gamma(1/2)} - \frac{1}{p}\\
0 &0 &\frac{1}{\Gamma(1/2)} &p^{-1}E_{\frac{1}{2},\frac{1}{2}}(p) - \frac{1}{p\Gamma(1/2)}\\
0 &0 &0 &E_{\frac{1}{2},\frac{1}{2}}(p)\\
\end{pmatrix}.
\end{gather*}

And taking into account the following properties of gamma-function and scalar Mittag-Leffler function
\begin{equation*}
\begin{aligned} 
&\Gamma\left( \frac{1}{2}\right) =\sqrt{\pi},  \Gamma\left( \frac{3}{2}\right)=\frac{1}{2}\Gamma\left( \frac{1}{2}\right)=\frac{\sqrt{\pi}}{2}, E_{\frac{1}{2},1}(z)=e^{z^2}\erfc(-z), \\
&E_{\frac{1}{2},\frac{1}{2}}\left(z\right)=z E_{\frac{1}{2},1}(z)+\frac{1}{\Gamma(1/2)} = z e^{z^2}\erfc(-z) + \frac{1}{\sqrt{\pi}}, 
\end{aligned}
\end{equation*}
we arrive at the following explicit expressions
\begin{equation}\label{H1}
\begin{gathered}
E_{\frac{1}{2},1}(B) = 
\begin{pmatrix} 
1 &2\sqrt{\frac{1}{\pi}} &1 &\frac{e^{p^2}}{p^3}\erfc(-p)- \frac{1}{p^3} - \frac{2}{p^2}\sqrt{\frac{1}{\pi}} - \frac{1}{p}\\
0 &1 &2\sqrt{\frac{1}{\pi}} &\frac{e^{p^2}}{p^2}\erfc(-p)- \frac{1}{p^2} - \frac{2}{p}\sqrt{\frac{1}{\pi}}\\
0 &0 &1 &\frac{e^{p^2}}{p}\erfc(-p)- \frac{1}{p}\\
0 &0 &0 &e^{p^2}\erfc(-p)\\
\end{pmatrix},
\end{gathered}
\end{equation}

\begin{equation}\label{H2}
\begin{gathered}
E_{\frac{1}{2},\frac{1}{2}}\left(B\right) = 
\begin{pmatrix} 
\frac{1}{\sqrt{\pi}} &1 &\frac{2}{\sqrt{\pi}} &\frac{e^{p^2}}{p^2}\erfc(-p) - \frac{1}{p^2} - \frac{2}{p\sqrt{\pi}}\\
0 &\frac{1}{\sqrt{\pi}} &1 &\frac{e^{p^2}}{p}\erfc(-p) - \frac{1}{p}\\
0 &0 &\frac{1}{\sqrt{\pi}} &e^{p^2} \erfc(-p)\\
0 &0 &0 &pe^{p^2}\erfc(-p)+\frac{1}{\sqrt{\pi}}\\
\end{pmatrix}.
\end{gathered}
\end{equation}

Here $ \erfc$ stands for the complementary error function, an entire function defined by 
\[
\erfc(z)=\frac{2}{\sqrt{\pi}}\int_z^\infty e^{-t^2}dt.
\]

The following \textsc{Matlab} code evaluates the matrix Mittag-Leffler function of $B$ for the case $ a=b$, $ c=0$, with the help of \verb|mlfm.m| routine and compares it with the reference matrices \eqref{H1}, \eqref{H2}

\begin{lstlisting}
B = [0  1  0  0;
     0  0  1  0;
     0  0  0  1;
     0  0  0 -1];

% matrix Mittag-Leffler function of B for alpha=0.5, beta=1
E1 = mlfm(B,0.5,1);                                        

% reference matrix
H1 = [1 2/sqrt(pi) 1          -exp(1)*erfc(1)-2/sqrt(pi)+2;
      0 1          2/sqrt(pi) exp(1)*erfc(1)+2/sqrt(pi)-1;
      0 0          1          -exp(1)*erfc(1)+1;
      0 0          0          exp(1)*erfc(1)];             

% matrix Mittag-Leffler function of B for alpha=0.5, beta=0.5
E2 = mlfm(B,0.5,0.5);                                      

% reference matrix
H2 = [1/sqrt(pi) 1          2/sqrt(pi) exp(1)*erfc(1)-1+2/sqrt(pi);
      0          1/sqrt(pi) 1          -exp(1)*erfc(1)+1;
      0          0          1/sqrt(pi) exp(1)*erfc(1);
      0          0          0          -exp(1)*erfc(1)+1/sqrt(pi)];
      
tol = 1e-15; % tolerance
abs(E1-H1) < tol
abs(E2-H2) < tol
\end{lstlisting}
The above code produces the result
\begin{verbatim}
ans =

1     1     1     1
1     1     1     1
1     1     1     1
1     1     1     1


ans =

1     1     1     1
1     1     1     1
1     1     1     1
1     1     1     1
\end{verbatim}
This implies that the \verb|mlfm.m| routine evaluates the matrix Mittag-Leffler function with sufficiently high accuracy (absolute error is less than $ 10^{-15} $).


\end{document}